\documentclass{article}
\usepackage{geometry}
\usepackage[utf8]{inputenc}
\usepackage{amsmath}
\usepackage{bm}
\usepackage{graphicx}
\usepackage{authblk}
\usepackage[autostyle]{csquotes} 
\usepackage{cite}
\usepackage[dvipsnames]{xcolor}
\usepackage[hidelinks]{hyperref}
\hypersetup{
    colorlinks = true,
    linkcolor={Blue}, citecolor={Blue}}
\title{$3n+3^k$: New Perspective on Collatz Conjecture}
\author{Naouel Boulkaboul\footnote{Email: nwl\_qcd@yahoo.com}\\ \emph{Independent researcher, Algiers, Algeria}}
\date{October 2022}

\begin{document}

\maketitle
\begin{abstract}
 Collatz conjecture is generalized to $3n+3^k$ ($k\in N$). Operating as usual, every sequence seems to reach $3^k$ and end up in the loop $3^k, 4.3^k, 2.3^k,3^k$. The usual $3n+1$ conjecture is recovered for $k=0$. For $k>0$, we noticed the existence of a sequence of period 3, namely, $3^{k-1}, 2.3^k, 3^k$, alongside the cycle $4.3^k, 2.3^k,3^k$ encountered in the $3n+1 (k=0)$ sequence. A term formula of the $3n+3^k$ conjecture has been derived, and hence the \emph{total stopping time}.    
\end{abstract}
\section{Introduction}
Credited to the mathematician Lothar Collatz, the Collatz conundrum (see for instance Ref. \cite{lagarias20033x+} and references therein), which was brought forward in the 1930s, is one of the simplest yet unsolved conjectures in mathematics. \emph{``This is a really dangerous problem. People become obsessed with it and it really is impossible,”} as stated by Jeffrey Lagarias. Such a conjecture asks whether repeating two simple arithmetic operations will eventually reach 1. That is to say, let $N:=\{0, 1, 2...\}$ denote the natural numbers, so
that $N+1=\{1, 2, 3...\}$ are the positive integers. Now, pick an arbitrary positive integer $n\in N+1$ and apply the following operation on it:
If the number is even, divide it by two; while if the number is odd, triple it and add one. Let the process be denoted by $f(n)$. That is
\begin{equation}
f(n)=
 \begin{cases}
 \frac{n}{2} \qquad \text{if} \quad n\equiv 0 \qquad(\text{mod} \quad 2)\\
3n+1 \qquad \text{if}\quad n\equiv 1 \qquad(\text{mod}\quad 2)
\end{cases}
\end{equation}
The main task is to form a sequence by performing this operation repeatedly, and take the output at each step as the input at the next. Consequently, one may notice that such a process will eventually lead to the number 1, regardless of which positive integer is chosen initially. \\
If we set $f^0(n)=n$ and $f^l(n)=f(f^{l-1}(n))$ for $l\in N$. Then, the Collatz sequence for $n$ reads
\begin{equation}
C(n)=\{f^l(n)\}^{\infty}_{l=0},
\end{equation}
Strictly speaking, one may note that every Collatz sequence ends up in the loop $1, 4, 2, 1$. To disprove Collatz conjecture, one has to show that there exists some starting number which yields a sequence that does not include 1. Such a sequence would either enter a repeating cycle that excludes 1, or increase infinitely. No such sequence has been found so far \cite{barina2021convergence}.\\
In the present paper, we show that Collatz conjecture $3n+1$ may possibly be a part of a more generalized case of the form $3n+3^k$ with $k\in N$.
\section{$3n+3^k$ Conjecture}
We provide the following operation, pick an arbitrary positive integer $n\in N+1$, if the number is even, divide it by two; otherwise if the number is odd, triple it and add $3^k$ with $k\in N$.
\begin{equation}
g_k(n)=
 \begin{cases}
 \frac{n}{2} \qquad \text{if} \quad n \equiv 0 \qquad(\text{mod} \quad 2)\\
3n+3^k \qquad \text{if}\quad n \equiv 1 \qquad(\text{mod} \quad 2) \qquad k\in N
\end{cases}
\end{equation}
Repeating the two operations will eventually lead to $3^k$ ($k\in N$). Hence, one may argue that the original Collatz conjecture $3n+1$ is nothing more than $3n+3^k$ conjecture with $k=0$, i. e. $f(n)=g_0(n)$. Apparently, each positive integer, upon repetitive application of $g_k(n)$, will end in a repeating sequence of the form $3^k, 4.3^k, 2.3^k, 3^k$. In the $3n+1$ (or equivalently $3n+3^k$ with $k=0$) sequence, the only cycle of period $3$ is known to be the sequence $4.3^k, 2.3^k, 3^k$ \cite{zarnowski2001generalized}. For $k>0$ though, along with the aforementioned cycle,  one may distinguish another sequence of period $3$, that is $3^{k-1}, 2.3^k, 3^k$. Although such a generalized conjecture seems to be true, it has yet to be numerically verified for at least the largest set of numbers that has been checked in the case of $3n+1$ so far, i. e. $n\leq 2^{100000} - 1$ \cite{ren2018collatz}, even though a general proof of it is still lacking. Furthermore, up to this point it is not clear whether there exists a certain value of $k$ for or beyond which the conjecture cease to be hold.\\
\section{Term formula for $3n+3^k$ sequence}
Let $C^l_k(n)$ include the first $l$ terms of the $3n+3^k$  sequence for $n$. We set $m$ as the number of odd terms in $C^l_k(n)$, $d_i$ as the number of consecutive even terms immediately following the $i$th odd term and $d_0$ as the number of even terms preceding the first odd term. It follows that the next term in the $3n+3^k$ sequence for $n$ is \footnote{Note that this formula is different from the one derived by L E. Garner \cite{garner1981collatz}}.
\begin{equation}
\label{next term}
g^l_k(n)=\frac{3^m}{2^{l-m}}n+3^k\epsilon\sum_{j=1}^{m}\frac{3^{m-j}}{2^{\sum_{i=j}^{m}d_i}}    
\end{equation}
where $\epsilon=0$ if $m=0$ and $\epsilon=1$ if $m\neq0$, note that $l-m=d_0+d_1+...+d_m$. It is worth noting that, for the original $3n+1$ Collatz sequence, the above relation reduces to
\begin{equation}
g^l_0(n)=\frac{3^m}{2^{l-m}}n+\epsilon\sum_{j=1}^{m}\frac{3^{m-j}}{2^{\sum_{i=j}^{m}d_i}}    
\end{equation}
$3n+3^k$ sequences of integers ranging from $1$ to $17$ for $k=0, 1, 2, 3$ and $4$ are listed in Table. \ref{table}. In the $k$-sequences chosen, it is evident that, after consecutive application of $g^l_k$, the aforementioned integers eventually reach $3^k$.
\section{Total stopping time}
\begin{figure*}
\begin{centering}
\includegraphics[scale=0.4]{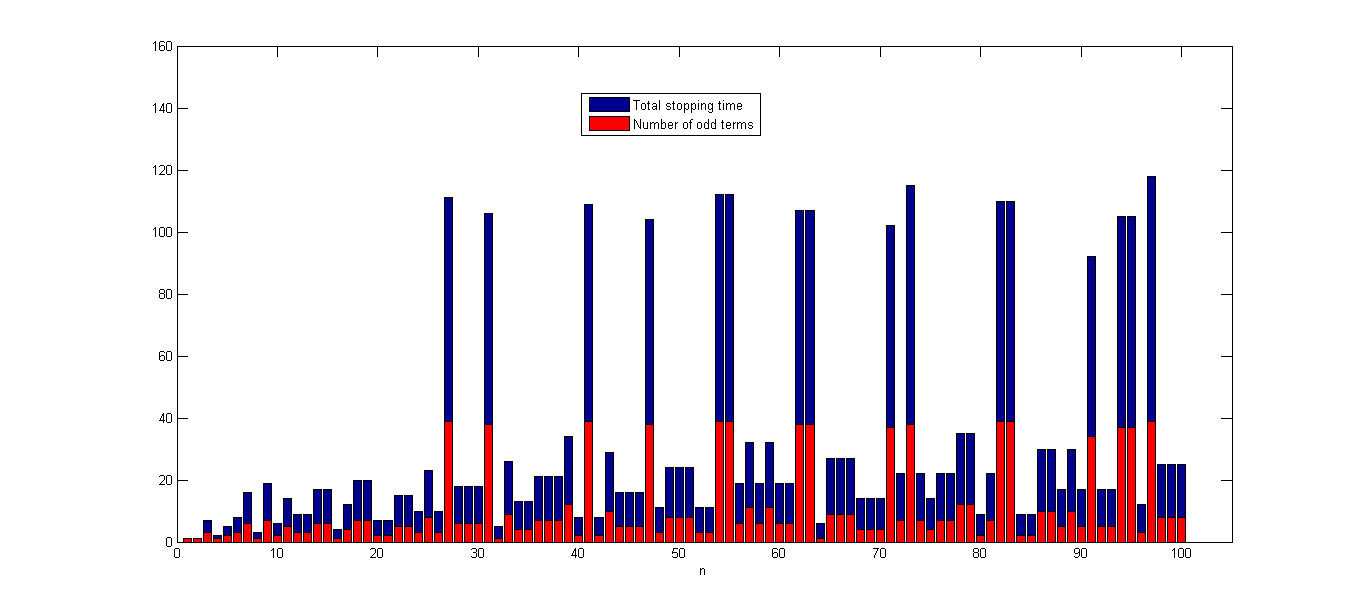}\\
\includegraphics[scale=0.4]{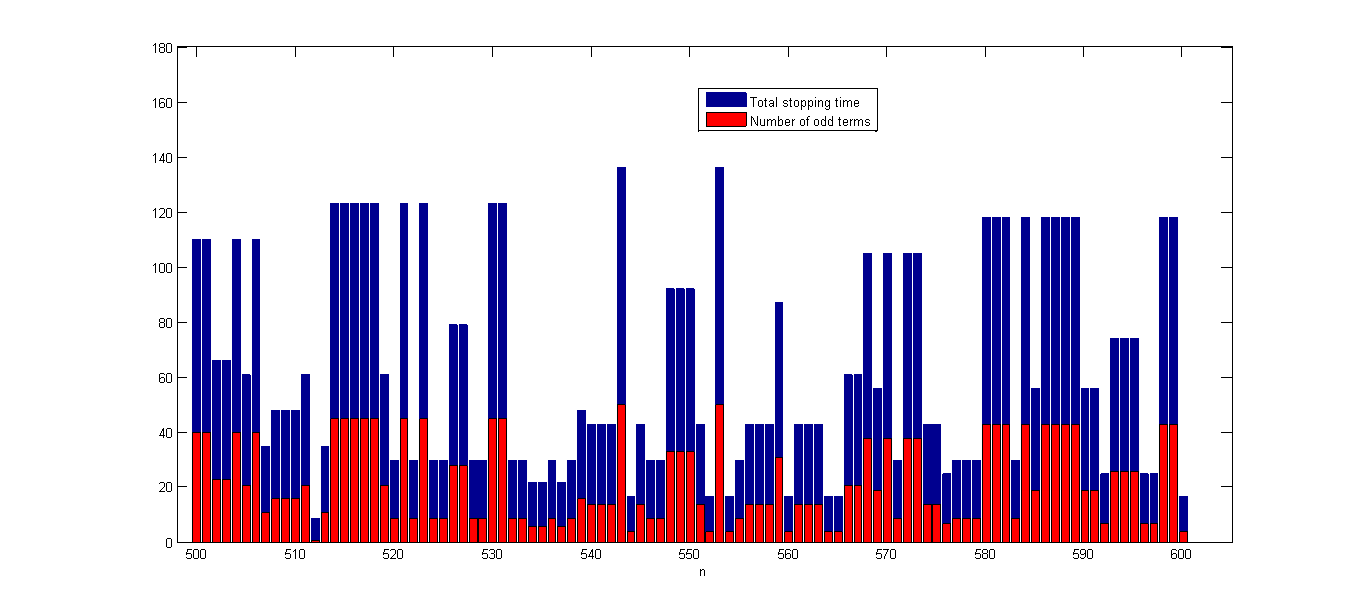}\\
\includegraphics[scale=0.4]{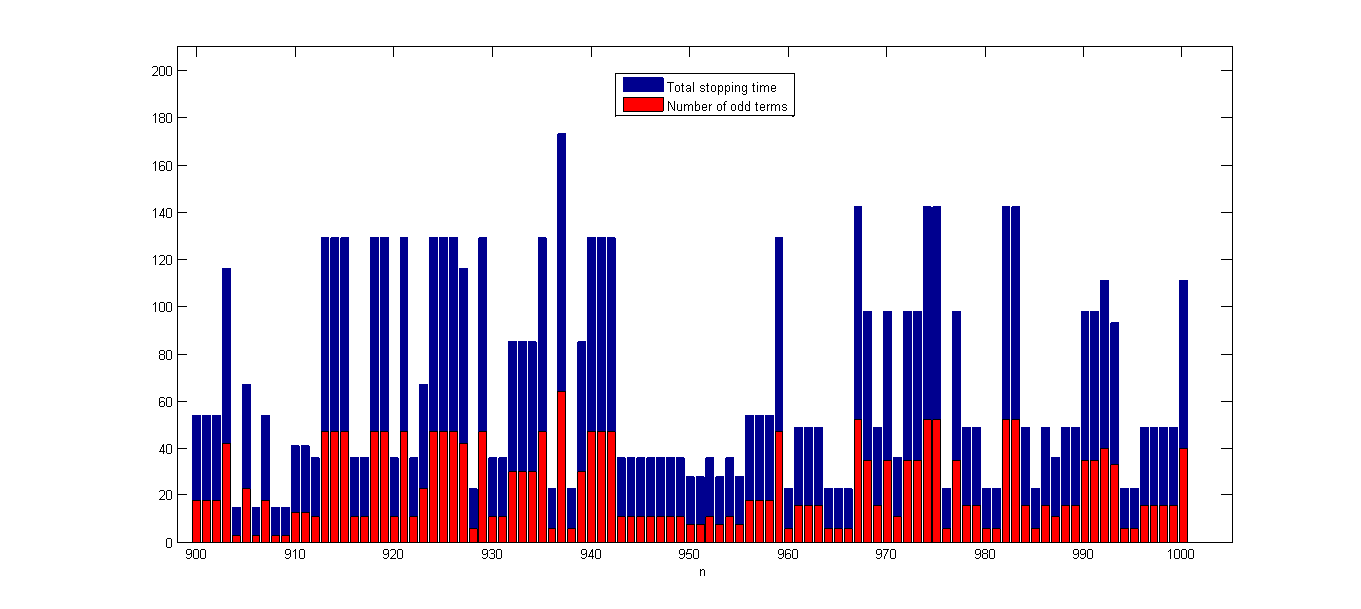}\\
\end{centering}
\caption{\label{fig:1}Behavior of the \emph{total stopping time} and the number of odd terms in the set $C^{t+1}_1(n)$ of the $3n+1$ sequence for $n\in[1,100]$ (top), $n\in[500,600]$ (middle) and $n\in[900,1000]$ (bottom), the blue bars represent the \emph{total stopping time} while the red bars depict the number of odd terms in the sequence.}
\end{figure*}

\begin{figure*}
\begin{centering}
\includegraphics[scale=0.4]{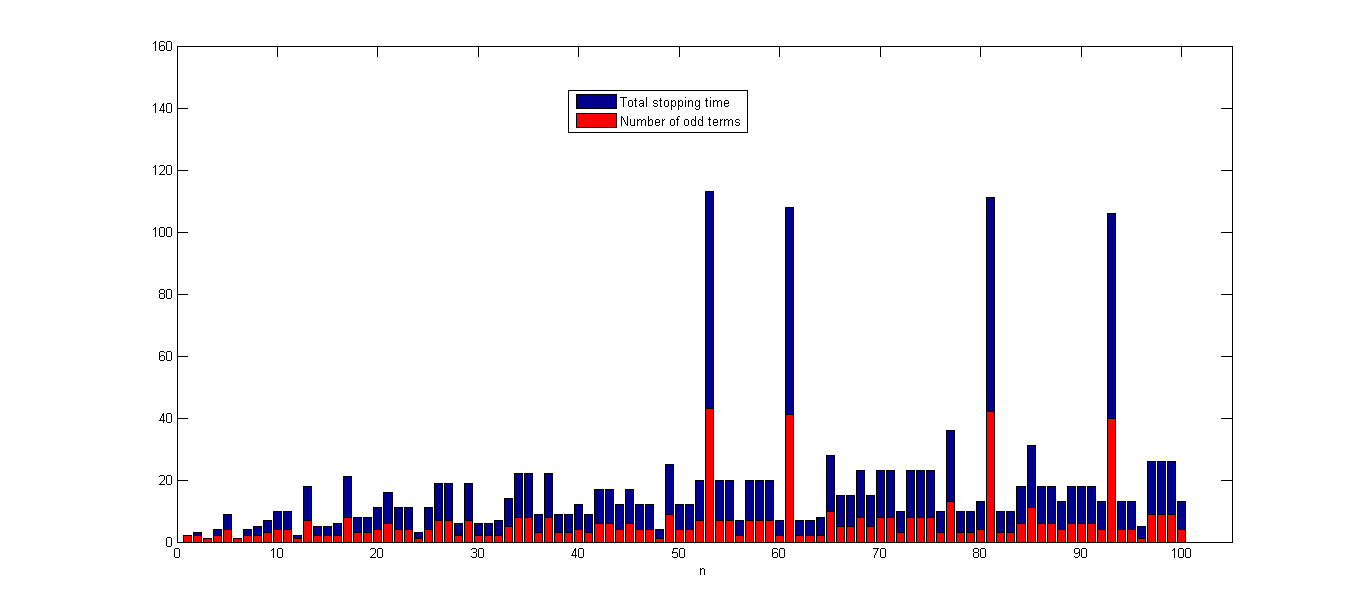}\\
\includegraphics[scale=0.4]{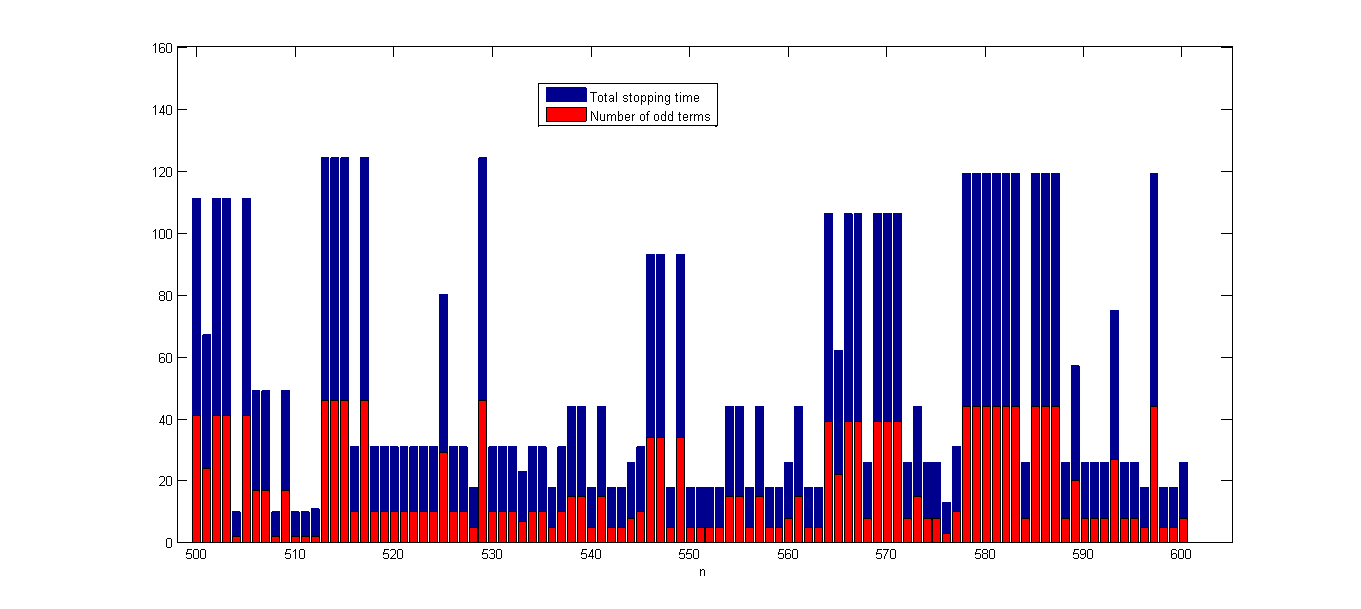}\\
\includegraphics[scale=0.4]{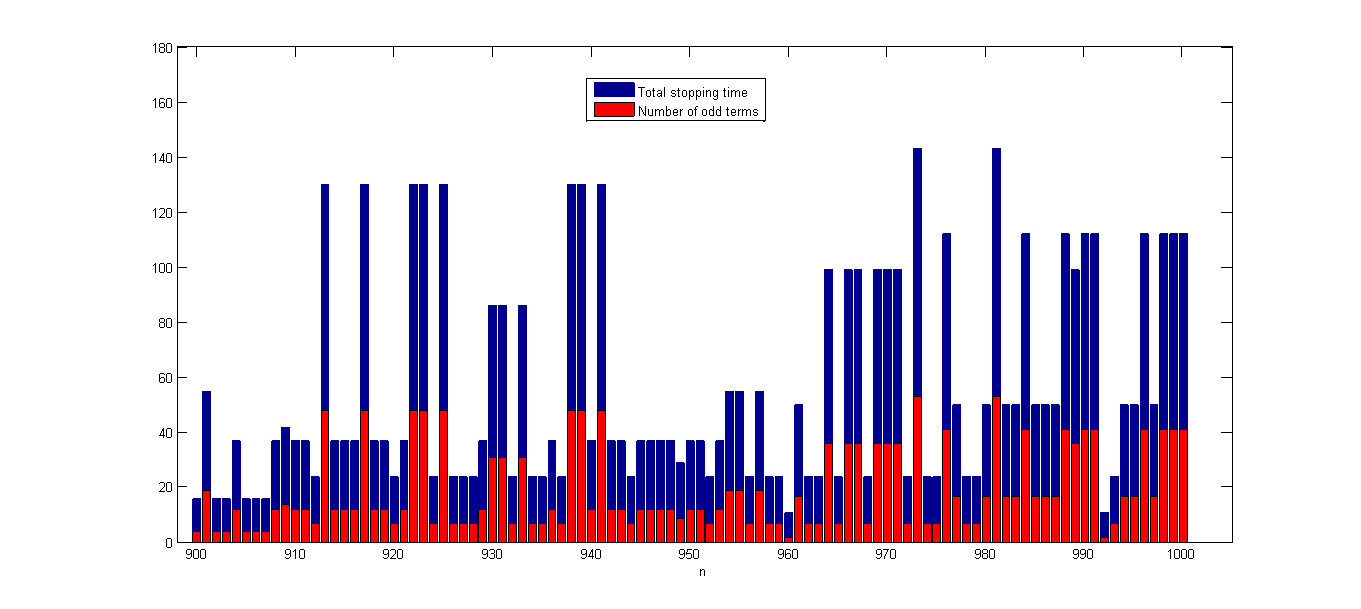}\\
\caption{\label{fig:2}Behavior of the \emph{total stopping time} and the number of odd terms in the set $C^{t+1}_3(n)$ of the $3n+3$ sequence for $n\in[1,100]$ (top), $n\in[500,600]$ (middle) and $n\in[900,1000]$ (bottom), the blue bars represent the \emph{total stopping time} while the red bars depict the number of odd terms in the sequence.}
\end{centering}
\end{figure*}

\begin{figure*}
\begin{centering}
\includegraphics[scale=0.4]{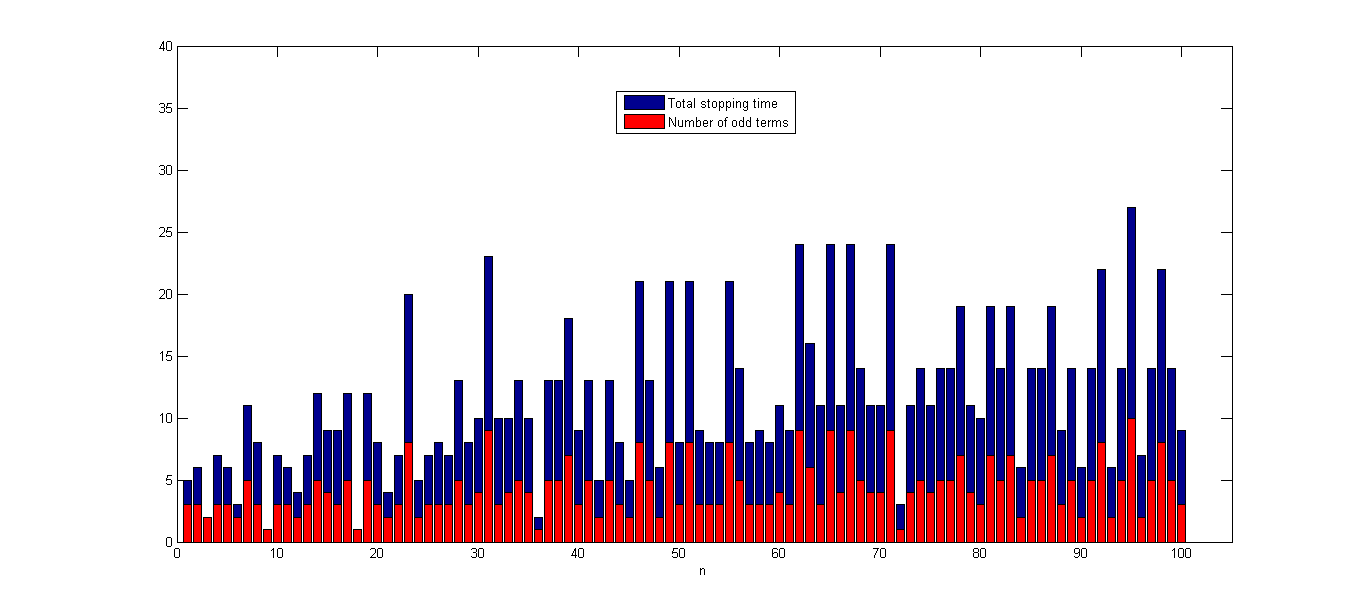}\\
\includegraphics[scale=0.4]{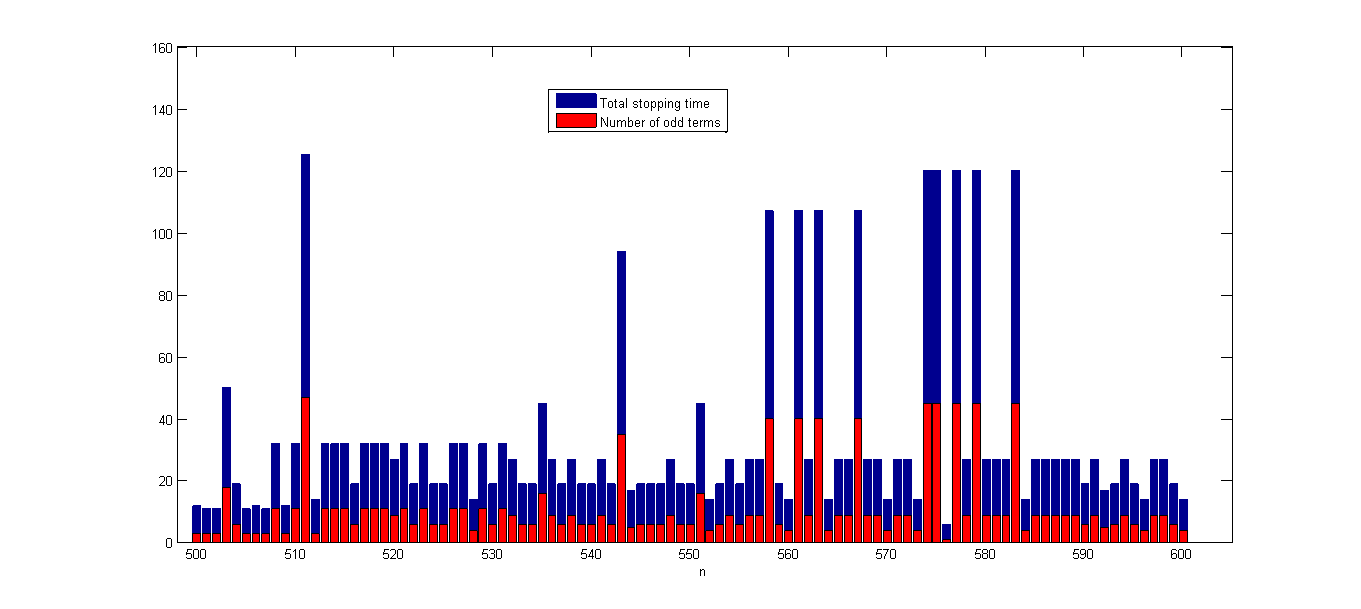}\\
\includegraphics[scale=0.4]{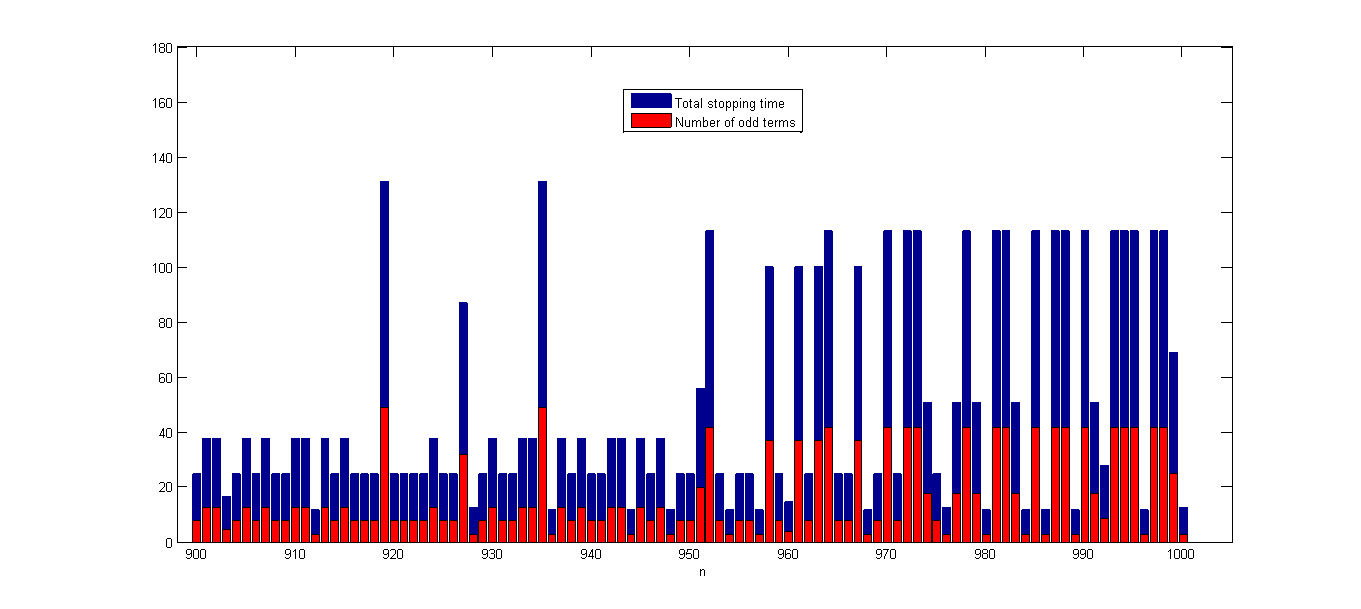}\\
\caption{\label{fig:3}Behavior of the \emph{total stopping time} and the number of odd terms in the set $C^{t+1}_9(n)$ of the $3n+9$ sequence for $n\in[1,100]$ (top), $n\in[500,600]$ (middle) and $n\in[900,1000]$ (bottom), the blue bars represent the \emph{total stopping time} while the red bars depict the number of odd terms in the sequence.}
\end{centering}
\end{figure*}

We call the \emph{total stopping time} of $n$ the smallest $t$ such that $n_t=3^k$, with $n_t$ being the value of $g_k$ applied to $n$ recursively $t$ times. It is worth mentioning that $t$ is nothing but the number of integers in the $3n+3^k$ sequence just preceding the $3^k$ term, i. e. $n_t=g^t_k(n)$. Thus, using the relation \eqref{next term}, the \emph{total stopping time} can be written as 
\begin{equation}
\label{total stoping time}
t=\text{log}_2\bigg[\frac{2^m3^m}{3^k(1-\epsilon\sum_{j=1}^{m}\frac{3^{m-j}}{2^{\sum_{i=j}^{m}d_i}})}  n\bigg]
\end{equation}
If $t$ does not exist we say that the \emph{total stopping time} is infinite. It follows that, one has to prove that all positive integers yield a finite \emph{total stopping time} in order to prove the $3n+3^k$ conjecture. In an endeavor to prove the $3n+1$ conjecture, mathematicians used to inspect what it is called the \emph{stopping time}, namely the least positive $l$ for which $f^l(n)=g^l_0(n)<n$ \cite{korec1994density, crooks2022collatz,idowu2015novel,ruggiero2019relationship}. If they can prove that
all positive integers have a finite \emph{stopping time}, they can prove by induction that the Collatz conjecture is true. It is nothing short of reasonable since for $n>1$, $f^l(n)=g^l_0(n)=1$ cannot occur without the occurrence of some $f^l(n)=g^l_0(n)<n$  \cite{lagarias2010ultimate}. Indeed, in the 1970s, it has been shown that almost all Collatz sequences eventually reach a number that is smaller than where it started \cite{terras1976stopping}. In 2019, Terence Tao proved that for almost all numbers the Collatz sequence of $n$ leads to a lower value, showing that the Collatz conjecture holds true for almost all numbers \cite{tao2019almost}.  But the question to be asked is whether it is the right (or sufficient) way to prove the Collatz conjecture. 
 As obvious as it seems, if the $3n+3^k$ conjecture turns out to be true for all $n\in N$+1, then the statement that, for each $n\in N+1$  there exists $l\in N+1$ such that $g^l_k(n)<n$, does not hold for $n<3^k$. This can be easily seen in the $3n+27$ sequence for $n=5, 9, 13$ (see for instance Table. \ref{table}). Alternatively, we must focus on showing that each $n\in N+1$ has a finite \emph{total stopping time} and hence proving that the $3n+3^k$ sequence inclines toward $3^k$ no matter what integer we start with. That is to say, one has to test the validity of the following statement: \emph{``for each $n\in N+1$, there exists $t\in N+1$ such that $g^t_k(n)=3^k$."}, without making use of the \emph{stopping time}.\\
 Now, suppose that two integers $n_1$ and $n_2$ have the same \emph{total stopping time} $t$ within the sequence $3n+3^k$, then from Eq. \eqref{total stoping time} we have 
\begin{equation}
\label{total stoping time2}
n_2=\frac{2^{m_1}3^{m_1}}{2^{m_2}3^{m_2}}\frac{1-\epsilon\sum_{j=1}^{m_2}\frac{3^{m_2-j}}{2^{\sum_{i=j}^{m_2}d_{2,i}}}}{1-\epsilon\sum_{j=1}^{m_1}\frac{3^{m_1-j}}{2^{\sum_{i=j}^{m_1}d_{1,i}}}}n_1
\end{equation}
where $m_1(m_2)$ and $d_{1,i}(d_{2,i})$ are, respectively, the number of odd terms and the number of consecutive even terms immediately following the $i$th odd term in $C^t_k(n_1)$ ($C^t_k(n_2)$), namely the set of $t$ terms of the $3n+3^k$ sequence for $n_1(n_2)$.\\
It is obvious that for $m_1=m_2=0$, we have $n_2=n_1$. Hence, no two integers, for which the corresponding sequences contain even terms only, have the same \emph{total stopping time}.  Sequences that have the same \emph{total stopping time} and the same cycle of period $3$ have the same number of odd terms, a feature that can be figured out by checking Figs. \ref{fig:1}$-$\ref{fig:3} which depict the \emph{total stopping time} and the number of odd terms  in the set $C^{t+1}_k(n) (k=0, 1, 2)$ for different values of $n$. In the $3n+9$ sequence for example (upper plot in Fig. \ref{fig:3}), there exist two cycles of period $3$: $36, 18, 9$ and $3, 18, 9$, the integers $32, 33$ and $35$ have the same \emph{total stopping time}, i.e. $t=10$, but only $33$ and $35$ share the same number of odd terms $m+1=4$, which is due to the fact that both of them possesses the same cycle of period $3$, i.e. $36, 18, 9$, whereas $32$ has the cycle $3, 18, 9$.\\
Remarkably, It is worth noting that sequences that possess the same number of odd terms in the set $C^{t+1}_k(n)$ do not necessarily yield the same \emph{total stopping time}, even if they share the same cycle of period 3, take the example of 1 and 2 in the $3n+9$ sequence, both of them have the same number of odd terms in the set $C^{t+1}_9(n)$, i.e. $m+1=3$, and the same cycle $3, 18, 9$ (check Table. \ref{table}) meanwhile they have different \emph{total stopping time}, $t=5$ and $t=6$, respectively.
\begin{table}
\begin{centering}
\begin{tabular}{lccccccccccccccccc}
   
\hline
\hline 
&&&&&&&&$3n+1$&&&&&&&&\\
\hline
$n$ & 1 & 2 & 3 & 4 & 5 & 6 & 7 & 8 & 9 & 10&11 &12 &13 &14 &15 &16 &17 \\
&4	&1	&10	&2	&16	&3	&22	&4	&28	&5	&34	&6	&40	&7	&46	&8	&52\\
&2	&4	&5	&1	&8	&10	&11	&2	&14	&16	&17	&3	&20	&22	&23	&4	&26\\
&1	&2	&16	&4	&4	&5	&34	&1	&7	&8	&52	&10	&10	&11	&70	&2	&13\\
&	&1	&8	&2	&2	&16	&17	&4	&22	&4	&26	&5	&5	&34	&35	&1	&40\\
&	&	&4	&1	&1	&8	&52	&2	&11	&2	&13	&16	&16	&17	&106 &4	&20\\
&	&	&2	&	&4	&4	&26	&1	&34	&1	&40	&8	&8	&52	&53	&2	&10\\
&	&	&1	&	&2	&2	&13	&	&17	&4	&20	&4	&4	&26	&160	&1	&5\\
&	&	&4	&	&1	&1	&40	&	&52	&2	&10	&2	&2	&13	&80	&	&16\\
&	&	&2	&	&	&4	&20	&	&26	&1	&5	&1	&1	&40	&40	&	&8\\
&	&	&1	&	&	&2	&10	&	&13	&	&16	&4	&4	&20	&20	&	&4\\
&	&	&	&	&	&1	&5	&	&40	&	&8	&2	&2	&10	&10	&	&2\\
&	&	&	&	&	&	&16	&	&20	&	&4	&1	&1	&5	&5	&	&1\\
&	&	&	&	&	&	&8	&	&10	&	&2	&	&	&16	&16	&	&4\\
&	&	&	&	&	&	&4	&	&5	&	&1	&	&	&8	&8	&	&2\\
&	&	&	&	&	&	&2	&	&16	&	&4	&	&	&4	&4	&	&1\\
&	&	&	&	&	&	&1	&	&8	&	&2	&	&	&2	&2	&	&\\
&	&	&	&	&	&	&4	&	&4	&	&1	&	&	&1	&1	&	&\\
&	&	&	&	&	&	&2	&	&2	&	&	&	&	&4	&4	&	&\\
&	&	&	&	&	&	&1	&	&1	&	&	&	&	&2	&2	&	&\\
&	&	&	&	&	&	&	&	&4	&	&	&	&	&1	&1	&	&\\
&	&	&	&	&	&	&	&	&2	&	&	&	&	&	&	&	&\\
&	&	&	&	&	&	&	&	&1	&	&	&	&	&	&	&	&\\
\hline

\hline 
&&&&&&&&$3n+3$&&&&&&&&\\
\hline
$n$ & 1 & 2 & 3 & 4 & 5 & 6 & 7 & 8 & 9 & 10&11 &12 &13 &14 &15 &16 &17 \\
&6	&1	&12	&2	&18	&3	&24	&4	&30	&5	&36	&6	&42	&7	&48	&8	&54\\
&3	&6	&6	&1	&9	&12	&12	&2	&15	&18	&18	&3	&21	&24	&24	&4	&27\\
&12	&3	&3	&6	&30	&6	&6	&1	&48	&9	&9	&12	&66	&12	&12	&2	&84\\
&6	&12	&	&3	&15	&3	&3	&6	&24	&30	&30	&6	&33	&6	&6	&1	&42\\
&3	&6	&	&12	&48	&	&12	&3	&12	&15	&15	&3	&102	&3	&3	&6	&21\\
&	&3	&	&6	&24	&	&6	&12	&6	&48	&48	&	&51	&12	&12	&3	&66\\
&	&	&	&3	&12	&	&3	&6	&3	&24	&24	&	&156	&6	&6	&12	&33\\
&	&	&	&	&6	&	&	&3	&12	&12	&12	&	&78	&3	&3	&6	&102\\
&	&	&	&	&3	&	&	&	&6	&6	&6	&	&39	&	&	&3	&51\\
&	&	&	&	&12	&	&	&	&3	&3	&3	&	&120	&	&	&	&156\\
&	&	&	&	&6	&	&	&	&	&12	&12	&	&60	&	&	&	&78\\
&	&	&	&	&3	&	&	&	&	&6	&6	&	&30	&	&	&	&39\\
&	&	&	&	&	&	&	&	&	&3	&3	&	&15	&	&	&	&120\\
&	&	&	&	&	&	&	&	&	&	&	&	&48	&	&	&	&60\\
&	&	&	&	&	&	&	&	&	&	&	&	&24	&	&	&	&30\\
&	&	&	&	&	&	&	&	&	&	&	&	&12	&	&	&	&15\\
&	&	&	&	&	&	&	&	&	&	&	&	&6	&	&	&	&48\\
&	&	&	&	&	&	&	&	&	&	&	&	&3	&	&	&	&24\\
&	&	&	&	&	&	&	&	&	&	&	&	&12	&	&	&	&12\\
&	&	&	&	&	&	&	&	&	&	&	&	&6	&	&	&	&6\\
&	&	&	&	&	&	&	&	&	&	&	&	&3	&	&	&	&3\\
&	&	&	&	&	&	&	&	&	&	&	&	&	&	&	&	&12\\
&	&	&	&	&	&	&	&	&	&	&	&	&	&	&	&	&6\\
&	&	&	&	&	&	&	&	&	&	&	&	&	&	&	&	&3\\

\hline
\end{tabular}
\end{centering}
\end{table}

\begin{table}
\begin{centering}
\begin{tabular}{lccccccccccccccccc}
   
\hline
\hline 
&&&&&&&&$3n+9$&&&&&&&&\\
\hline
$n$ & 1 & 2 & 3 & 4 & 5 & 6 & 7 & 8 & 9 & 10&11 &12 &13 &14 &15 &16 &17 \\
&12	&1	&18	&2	&24	&3	&30	&4	&36	&5	&42	&6	&48	&7	&54	&8	&60\\
&6	&12	&9	&1	&12	&18	&15	&2	&18	&24	&21	&3	&24	&30	&27	&4	&30\\
&3	&6	&36	&12	&6	&9	&54	&1	&9	&12	&72	&18	&12	&15	&90	&2	&15\\
&18	&3	&18	&6	&3	&36	&27	&12	&36	&6	&36	&9	&6	&54	&45	&1	&54\\
&9	&18	&9	&3	&18	&18	&90	&6	&18	&3	&18	&36	&3	&27	&144 &12	&27\\
&36	&9	&	&18	&9	&9	&45	&3	&9	&18	&9	&18	&18	&90	&72	&6	&90\\
&18	&36	&	&9	&36	&	&144 &18	&	&9	&36	&9	&9	&45	&36	&3	&45\\
&9	&18	&	&36	&18	&	&72	&9	&	&36	&18	&	&36	&144	&18	&18	&144\\
&	&9	&	&18	&9	&	&36	&36	&	&18	&9	&	&18	&72	&9	&9	&72\\
&	&	&	&9	&	&	&18	&18	&	&9	&	&	&9	&36	&36	&36	&36\\
&	&	&	&	&	&	&9	&9	&	&	& &	&	&18	&18	&18	&18\\
&	&	&	&	&	&	&36	&	&	&	&	&	&	&9	&9	&9	&9\\
&	&	&	&	&	&	&18	&	&	&	&	&	&	&36	&	&	&36\\
&	&	&	& &	&	&9	&	&	&	&	&	&	&18	&	&	&18\\
&	&	&	&	&	&	&	&	&	&	&	&	&  &9	&	&	&9\\
\hline

\hline 
&&&&&&&&$3n+27$&&&&&&&&\\
\hline
$n$ & 1 & 2 & 3 & 4 & 5 & 6 & 7 & 8 & 9 & 10&11 &12 &13 &14 &15 &16 &17 \\
&30	&1	&36	&2	&42	&3	&48	&4	&54	&5	&60	&6	&66	&7	&72	&8	&78\\
&15	&30	&18	&1	&21	&36	&24	&2	&27	&42	&30	&3	&33	&48	&36	&4	&39\\
&72	&15	&9	&30	&90	&18	&12	&1	&108	&21	&15	&36	&126	&24	&18	&2	&144\\
&36	&72	&54	&15	&45	&9	&6	&30	&54	&90	&72	&18	&63	&12	&9	&1	&72\\
&18	&36	&27	&72	&162	&54	&3	&15	&27	&45	&36	&9	&216	&6	&54	&30	&36\\
&9	&18	&108	&36	&81	&27	&36	&72	&108	&162	&18	&54	&108	&3	&27	&15	&18\\
&54	&9	&54	&18	&270	&108	&18	&36	&54	&81	&9	&27	&54	&36	&108	&72	&9\\
&27	&54	&27	&9	&135	&54	&9	&18	&27	&270	&54	&108	&27	&18	&54	&36	&54\\
&108	&27	&	&54	&432	&27	&54	&9	&	&135	&27	&54	&108	&9	&27	&18	&27\\
&54	&108	&	&27	&216	&	&27	&54	&	&432	&108	&27	&54	&54	&	&9	&108\\
&27	&54	&	&108	&108	&	&108	&27	&	&216	&54	&	&27	&27	&	&54	&54\\
&	&27	&	&54	&54	&	&54	&108	&	&108	&27	&	&	&108	&	&27	&27\\
&	&	&	&27	&27	&	&27	&54	&	&54	&	&	&	&54	&	&108	&\\

&	&	&	&	&108	&	&	&27	&	&27	&	&	&	&27	&	&54	&\\
&	&	&	&	&54	&	&	&	&	&108	&	&	&	&	&	&27	&\\
&	&	&	&	&27	&	&	&	&	&54	&	&	&	&	&	&	&\\
&	&	&	&	&	&	&	&	&	&27	&	&	&	&	&	&	&\\
\hline
\end{tabular}
\end{centering}
\end{table}
\begin{table}
\label{collatz}
\begin{centering}
\begin{tabular}{lccccccccccccccccc}

\hline
\hline 
&&&&&&&&$3n+81$&&&&&&&&\\
\hline
$n$ & 1 & 2 & 3 & 4 & 5 & 6 & 7 & 8 & 9 & 10&11 &12 &13 &14 &15 &16 &17 \\
&84	&1	&90	&2	&96	&3	&102	&4	&108	&5	&114	&6	&120	&7	&126	&8	&132\\
&42	&84	&45	&1	&48	&90	&51	&2	&54	&96	&57	&3	&60	&102	&63	&4	&66\\
&21	&42	&216	&84	&24	&45	&234	&1	&27	&48	&252	&90	&30	&51	&270	&2	&33\\
&144	&21	&108	&42	&12	&216	&117	&84	&162	&24	&126	&45	&15	&234	&135	&1	&180\\
&72	&144	&54	&21	&6	&108	&432	&42	&81	&12	&63	&216	&126	&117	&486	&84	&90\\
&36	&72	&27	&144	&3	&54	&216	&21	&324	&6	&270	&108	&63	&432	&243	&42	&45\\
&18	&36	&162	&72	&90	&27	&108	&144	&162	&3	&135	&54	&270	&216	&810	&21	&216\\
&9	&18	&81	&36	&45	&162	&54	&72	&81	&90	&486	&27	&135	&108	&405	&144	&108\\
&108	&9	&324	&18	&216	&81	&27	&36	&	&45	&243	&162	&486	&54	&1296	&72	&54\\
&54	&108	&162	&9	&108	&324	&162	&18	&	&216	&810	&81	&243	&27	&648	&36	&27\\
&27	&54	&81	&108	&54	&162	&81	&9	&	&108	&405	&324	&810	&162	&324	&18	&162\\
&162	&27	&	&54	&27	&81	&324	&108	&	&54	&1296	&162	&405	&81	&162	&9	&81\\
&81	&162	&	&27	&162	&	&162	&54	&	&27	&648	&81	&1296	&324	&81	&108	&324\\
&324	&81	&	&162	&81	&	&81	&27	&	&162	&324	&	&648	&162	&324	&54	&162\\
&162	&324	&	&81	&324	&	&	&162	&	&81	&162	&	&324	&81	&162	&27	&81\\
&81	&162	&	&324	&162	&	&	&81	&	&324	&81	&	&162	&	&81	&162	&\\
&	&81	&	&162	&81	&	&	&324	&	&162	&324	&	&81	&	&	&81	&\\
&	&	&	&81	&	&	&	&162	&	&81	&162	&	&324	&	&	&324	&\\
&	&	&	&	&	&	&	&81	&	&	&81	&	&162	&	&	&162	&\\
&	&	&	&	&	&	&	&	&	&	&	&	&81	&	&	&81	&\\

\hline

\end{tabular}
\end{centering}
\caption{\label{table} $3n+3^k$ sequence with $k=0, 1, 2,3$ and $4$ for $n$ ranging from $1$ to $17$.}
\end{table}

\nocite{*}
\bibliographystyle{unsrt}
\bibliography{collatz}
\end{document}